\documentclass{article}
\usepackage{graphicx} 

\usepackage{amssymb,amsmath,enumerate}
\usepackage{xcolor}

\newtheorem{theorem}{Theorem}[section]

\newtheorem{conjecture}[theorem]{Conjecture}

\newcommand{\proofbox}{\hbox{\vbox{\hrule\hbox{k\vrule\phantom{\vrule height 8pt
                         width 5pt depth 0pt}\vrule}\hrule}\quad}}
                        
      \renewcommand{\epsilon}{\varepsilon}

\title{\bf 
Embedding 
Nearly Spanning Trees } 
\author{
Bruce Reed\footnote{Mathematical Institute, Academica Sinica, Taiwan.  (\texttt{bruce.al.reed@gmail.com}).} 
\quad
Maya Stein\footnote{Department of Mathematical Engineering and Center for Mathematical Modeling (CNRS IRL2807), University of Chile. Supported by FONDECYT Regular Grant 1221905,  by ANID Basal Grant CMM FB210005,  and by MSCA-RISE-2020-101007705 project {\it RandNET}. (\texttt{mstein@dim.uchile.cl}).}
}
\date{}
\begin{document}

\maketitle

    \begin{center}
        {\it 
    Dedicated to the memory of Vera T. S\'os. 
    }
    \end{center}

 \begin{abstract} 
 The Erd\H os-S\'os Conjecture states that every graph with average degree exceeding $k-1$ contains every tree with $k$ 
edges as a subgraph.   We prove that there are $\delta>0$ and $k_0\in\mathbb N$ 
such that the conjecture holds for every tree $T$ with $k \ge k_0$ edges  and every graph $G$ with $|V(G)| \le (1+\delta)|V(T)|$.  
\end{abstract}

 \section{Introduction}

One of the best known conjectures in extremal graph theory is the
 Erd\H os-S\'os  Conjecture (see~\cite{Erdos64}). 
 \begin{conjecture}[Erd\H os-S\'os  Conjecture]
    \label{conjES}
    Every graph $G$  with average degree  $d(G)>k-1$ contains every tree $T$  with $k$ 
edges as a subgraph. 
\end{conjecture}
The conjecture clearly holds for  stars  and is known to hold for
 paths~\cite{ErdGall59}. It also  holds for large trees  whose maximum degree is  bounded~\cite{PokHyper}: this result relies on earlier results for dense host graphs~\cite{BPS3, rozhon}.
  Further,
Conjecture \ref{conjES}   holds for host graphs that are  bipartite~\cite{kalai-bip} or have no $4$-cycles~\cite{sacwoz}, and it holds whenever $k\ge |V(G)| -c$,
where $c$ is any given constant and $k$ is sufficiently large depending on $c$~\cite{goerlich2016}. For more background see~\cite{tree-survey}.

We prove the high density case of Conjecture~\ref{conjES} for large $k$, with no restrictions at all on the host graph or on the tree. Our result reads as follows:

\begin{theorem}
    \label{smallGtheorem1}
    There are $k_0\in\mathbb N$ and  $\delta>0$ 
    such that  for 
    all $k \ge k_0$ every graph~$G$ with $|V(G)| \le (1+\delta)k$ and  with average degree  $d(G)>k-1$ contains every tree $T$  with $k$ 
edges as a subgraph. 
\end{theorem}

In our proof of this theorem, we take $\delta$ to be $10^{-10}$, and we make no effort to optimise this number.
Our proof of Theorem \ref{smallGtheorem1} relies crucially on an earlier result by the authors. This result, shown in~\cite{RS23a, RS23b}, is the (surprisingly hard to prove) spanning tree case of a more general conjecture of Havet, Wood and the authors from~\cite{HRSW}, which states that every graph of minimum degree at least $2k/3$ and maximum degree at least $k$ contains all trees with $k$ edges as subgraphs.

\begin{theorem}
\label{RS}$\!\!${\rm\bf\cite{RS23a, RS23b}}
There is an $m_0 \in\mathbb N $ such that for every $m \ge m_0$ every graph on $m + 1$ vertices that has minimum degree at least $\lfloor 2m/3 \rfloor$ and a  vertex of degree~$m$ contains every tree T with $m$ edges as a subgraph.
\end{theorem}

Let us quickly give some insight into the main ideas of our proof of Theorem~\ref{smallGtheorem1}. Since $G$ has high average degree, but no subgraph fulfilling the conditions of Theorem~\ref{RS}, we are able to find a relatively large set~$H\subseteq V(G)$ having degree   at least $k$. If $T$ has many leaves, we  use $H$ to embed parents of leaves while embedding the rest of the tree almost greedily, leaving it to the end to embed leaves from $H$. 
  
  The more difficult case is when $T$ has few leaves. Then $T$ has many vertices of degree~$2$, which  we embed at the end. The problem is that we have to choose the images of the neighbours of these vertices very carefully,  as we will need that at the very end of the embedding process, they are  adjacent to most of the yet unused vertices.  In order to achieve this, we  embed most  of the tree into a randomly selected set. More precisely, we embed most of the tree in an ordered way into a randomly selected permutation of most of $V(G)$. 
  For this approach to work, it is crucial that we have used up all vertices of $G$ having rather low degree in a different way in the beginning of the embedding process.

\section{The proof of Theorem \ref{smallGtheorem1}}

Set $\delta:=10^{-10}$. Let $m_0$ be given by  Theorem~\ref{RS} and  set 
$k_0:=\max \{m_0,\delta^{-2}\}.$
Let~$T$ be a tree with $k\ge k_0$ edges, and let 
$G$ be a graph with $n:=|V(G)| \le (1+\delta)k$ and $d(G)>k-1$. We can assume that $G$ is minimal with these properties. In particular, $G$ has no vertex $v$ of degree less than $\frac k2$, as deleting such a vertex would lead to a smaller graph of at least the same average degree. 
So,
     \begin{equation}\label{degL'}
       \delta (G)\ge \frac k2,
    \end{equation}
where $\delta(G)$, as usual, denotes the minimum degree of $G$.
Set $a:=n-k$,   $S\subseteq V(G)$ to be  the set  of all 
vertices of degree at most $\frac{2k}{3}+a$ (the letter $S$ stands for Small),
and  $b=|S|$. 

%

If $d(v)\ge k+b$ for some $v\in V(G)$, then $G_v:=G\big[N[v]-S\big]$  (the subgraph of $G$ induced by $v$ and all its neighbours outside $S$)
has at least $k+1$ vertices, each of which has at most $n-|V(G_v)| \le a -1$ neighbours outside~$G_v$. So the minimum degree of $G_v$ is at least $\frac{2k}{3}+a-(a-1)> \frac{2k}{3}$. Note that $v$ has degree $|V(G_v)|-1$ in  $G_v$, and hence we may apply
Theorem~\ref{RS} to find that $T\subseteq G_v\subseteq G$, and we are done. Therefore, we assume from now on that 
\begin{equation}\label{maxG}
    \Delta(G)< k+b.
\end{equation}
Let $H$ be the set of all vertices of $G$  having degree 
  at least $k$
  (the letter $H$ stands for High). Then
\begin{equation}\label{H}
    |H|> \frac{k}{6},
\end{equation}
as otherwise, \eqref{maxG} together with our choice of $\delta$,
ensures that 
\begin{align*} 
\sum_{v\in V(G)}d(v) \le \frac{k}{6} (k+b-1)+ (n-\frac k6-b)(k -1) +b(\frac{2k}3+a) \le n(k-1),
\end{align*}
a contradiction since $d(G)>k-1$.

Observe that the number of  vertices of $G$ having  degree less than 
    $k-\sqrt{ak}$  is at most $2 \sqrt{ak} \le 2\sqrt{\delta}k$ (as otherwise there are more than $ak\ge \frac{an}2$ non-edges and thus $d(\overline G)>a=n-k$, a contradiction).  
    We let $S'$ be the set of the $\lceil 2\sqrt{\delta}k \rceil$ vertices of lowest degrees. Then for each $v\in V(G)\setminus S'$, we have $d(v)\ge k-\sqrt{ak}$, and in particular, setting $G':=G[V(G)\setminus S']$, we have
    \begin{equation}\label{degG'}
        \delta(G')\ge k-\sqrt{ak}-|S'|\ge k-4\sqrt{\delta}k.
    \end{equation}

   Note that for each two vertices $u,v\in V(G')$, we have that
        \begin{equation}\label{codegG'}
       |N(u)\cap N(v) \cap V(G')|\ge  k- 9\sqrt\delta k.
    \end{equation}

 According to whether the tree has many or few leaves, we will either use $S'$ for embedding leaves at the end of our embedding procedure, or fill $S'$
    as early as possible.
For this, we distinguish two cases.\\

     {\bf Case 1:} $T$ contains at least $10\sqrt\delta k$ leaves.   \medskip

    Choose a set $L$ of exactly $\lceil 10\sqrt\delta k\rceil$ 
    leaves, and let $P_1$ be the set of their parents (note that not all children of $P_1$ need to belong to $L$). 
    Root $T$ in an arbitrary vertex $r$ and construct a set $P_2$ as follows.
    We start by setting $P_2:=P_1\cup \{r\}$. Then, while there is a vertex   $p\in V(T)\setminus P_2$ having at least 
    two children in $P_2$, we add $p$ to $P_2$. If there is no such~$p$, we stop the process. 
    We note that $|P_2|\le 2|P_1| \le 2|L|\le  21\sqrt{\delta}k$.  
    
    Let $P_3$ be the set of all parents of vertices from $P_2$.
    We embed $T[P_2]$ greedily  into 
    $H$, which is possible by~\eqref{H} and  the fact that the vertices of $H$ have degree at least $k$ in $G$ (and thus degree greater than $|P_2|$ into~$H$). We  
then embed $T-P_2-L$  into $G'$, going through $T-L$  in a top down fashion, starting with the root $r$, which, as it belongs to $P_2$, is already  embedded. 
For each subsequent vertex~$v$ that is not yet embedded, we choose its image  arbitrarily among the available neighbours of the image~$p_v$ of the parent  of $v$, unless $v\in P_3$, in which case we embed $v$ in a vertex that is adjacent to both $p_v$ and the image of the unique neighbour of $v$ in $P_2$. 
 This is possible by~\eqref{codegG'}, because $|L|\ge 10\sqrt\delta k$, and since so far, we only used $G'$ for the embedding. 
 
 It only remains to embed the vertices of $L$. Note that the parents of these vertices were embedded in $H$, and the vertices of $H$ have degree at least $k$. Thus we can embed $L$ greedily into $G$.\\

    {\bf Case 2:} $T$ has  fewer than $10\sqrt\delta k$ leaves.  \medskip
    
    In this case,  $T$ has fewer than $10\sqrt\delta k$ vertices of degree at least $3$. So the set~$D_2$  of vertices of degree $2$ has size greater than $k-20\sqrt\delta k$.  
We embed the vertices of $D_1:=V(T)\setminus D_2$ arbitrarily into $G'$ (respecting adjacencies) which is possible by~\eqref{degG'}. Let  $\varphi$ denote this embedding and all future extensions of it.

 Let $\mathcal R$ be the set of all components of $T[D_2]$. Note that each such component is a path and that 
    \begin{equation}\label{P}
 |\mathcal R|\le |D_1|\le 20\sqrt\delta k.
    \end{equation}
%
Take a minimal subset $\mathcal R'_1$ of $\mathcal R$ such that $\bigcup \mathcal R'_1$ contains at least $100\sqrt\delta k$ vertices. Choose an arbitrary path $Q$ from $\mathcal R'_1$ and delete one of its edges, giving us two subpaths $Q_1, Q_2$ of $Q$, in a way that $\mathcal R_1:=(\mathcal R'_1\setminus\{Q\})\cup\{Q_1\}$ covers exactly $\lfloor 100\sqrt\delta k\rfloor$ vertices. Set $\mathcal R_2:=(\mathcal R\setminus \mathcal R'_1)\cup \{Q_2\}$.
 
 For each path  $R\in\mathcal R_1$, we proceed as follows. Say $R=x_1x_2\ldots x_m$. 
 Set $X:=\{x_{2+3i}:0\le i\le m/3 -1 \}$ and note that $|X|\ge \lfloor 100\sqrt\delta k\rfloor - 4|\mathcal R_1|\ge |S'|$. 
 We embed an arbitrary subset $X'\subseteq X$ of size $|S'|$  arbitrarily into~$S'$. Then we embed the vertices from $V(R)\setminus X'$ into $G'$, in any order. Note that at the moment of being embedded, each such vertex $v$ has at most two already embedded neighbours, at most one of which is embedded in $S'$. By~\eqref{degL'} and~\eqref{degG'}, 
 there are at least $k-4\sqrt\delta k- (\frac k2 +a)=\frac k2- a-4\sqrt\delta k$ vertices that are common neighbours of the images of the neighbours of $v$.
 Moreover, we have embedded at most $|D_1|+\lfloor 100\sqrt\delta k\rfloor\le 120\sqrt\delta k$ vertices so far. Thus, 
 there are at least $$\frac{k}{2}-a-4\sqrt\delta k-120\sqrt\delta k>
 0
 $$ vertices of $G'$ that can serve as an appropriate image for $v$. Hence we can embed $\bigcup \mathcal R_1$ as planned.
 
 Let $U\subseteq V(G)$ be the set of all vertices used so far for the embedding. Note that 
 \begin{equation}\label{U}
 |U|\le  120\sqrt\delta k \le \frac k{300}. 
 \end{equation}
 It remains to embed the paths from $\mathcal R_2$ into $G-U$.
For this, let us introduce some notation.
 Given a  permutation $\pi=(v_1,v_2,\ldots, v_{n'})$ of $V(G)\setminus U$, we set $V_\pi:=\{v_1, v_2, \ldots, v_{\lceil\frac{49}{50}k\rceil}\}$.
 Let $J_\pi$ be the set of all indices $i< \lceil\frac{49}{50}k\rceil$ such that $v_i$ is not adjacent to $v_{i+1}$. 
 Let $H_\pi$ be the set of all vertices in $H\setminus (U\cup V_\pi)$ having less than $\frac{a}{3}$ non-neighbours in 
    $G'\setminus(U\cup V_\pi)$.    
    
     We claim that there is a permutation $\pi$ of $V(G)\setminus U$ such that
    \begin{enumerate}[(A)]
    \item  \label{A} $|J_\pi|\le 30\sqrt\delta k$, and 
    \item  \label{B} $|H_\pi|\ge 16 \sqrt\delta k$.
     \end{enumerate}
   Assuming such a permutation $\pi$ exists, we can finish the embedding as we will explain now. Choose any~$H'_\pi\subseteq H_\pi$ of size exactly $\lceil 16 \sqrt\delta k\rceil$ which is possible by~\eqref{B}. 
   We start by successively embedding paths from $\mathcal R_2$ as follows until we have used all of $V_\pi$. 
   We use the paths from $\mathcal R_2$  in non-decreasing order of their length.
    We embed each $R=x_1x_2\ldots x_m\in\mathcal R_2$ vertex by vertex, avoiding $H'_\pi$. Say we are at vertex $x_j\in V(R)$ with $j\neq m$. If possible we embed $x_j$ in the vertex $v_i$ with lowest index $i$ that has not been used yet. Otherwise we can and do embed~$x_j$ in a neighbour of $v_i\in V(G')\setminus H'_\pi$. 
     Vertex $x_m$ has  two already embedded neighbours $x, x'$ neither of which is embedded in $S'$, and  we embed $x_m$ in a common neighbour of  $\varphi(x)$ and  $\varphi(x')$, avoiding $H'_\pi$.
    All of this is possible by~\eqref{codegG'}, and since at any point, we have used 
    at most $2|\mathcal R_2|+|J_\pi|\le 70\sqrt\delta k$ vertices outside $V_\pi$, where the inequality holds by~\eqref{P} and~\eqref{A}. 
    We stop once we have used all of $V_\pi$, and let  $R'$ be the remainder of the path we were currently embedding. 
    Let $\mathcal R_3$ consist of $R'$ and all remaining paths of $\mathcal R_2$. Observe that since $|V(G)\setminus(U\cup V_\pi)|\le \frac k{50}$ and because of the order in which we used the paths from $\mathcal R_2$,
    \begin{equation}\label{P3}
     |\mathcal R_3| \le \frac{|\mathcal R_2|}{50}+1\le \frac{\sqrt\delta k}2.
      \end{equation}
      where the second inequality follows from~\eqref{P}. 
      
      By~\eqref{U},  the current total amount of vertices of $G$ used for the embedding is at most $|U|+|V_\pi|+70\sqrt\delta k \le \frac {99}{100}k-1$. Therefore, $$\big|\bigcup\mathcal R_3\big|\ge \frac {k}{100}\ge 2|H'_\pi|+3|\mathcal R_3|$$ and thus, there are sufficiently many vertices on the paths from $\mathcal R_3$ such that we can 
  embed the paths from $\mathcal R_3$ as follows. For each path $x_1x_2\ldots x_m\in\mathcal R_3$, we successively embed all vertices $x_j$ with even index  $j\neq m$ into $H'_\pi$, as long as there still are unused vertices in $H'_\pi$. For each odd index $j\notin\{1, m-1, m\}$ having the property that $x_{j-1}$ and $x_{j+1}$ are embedded in $H'_\pi$, we add vertex~$x_j$ to a set $W$, which is to be embedded at the very end.   
Observe that by construction, and by~\eqref{P3}, $$|W|\ge |H'_\pi|-|\mathcal R_3|\ge  15\sqrt{\delta}k.$$
We now use~\eqref{codegG'} to embed all remaining vertices from $V(\bigcup \mathcal R_3)\setminus W$ into~$G'$. 

 Finally, we embed $W$. By construction,   each vertex of $W$ is an $x_j$ from some path of $\mathcal R_3$, with $x_{j-1}, x_{j+1}$  embedded in vertices $u,v\in H'_\pi\subseteq H_\pi$. By definition of $H_\pi$,   at most $\frac{2a}{3}$ vertices in $V(G)\setminus (U\cup V_\pi)$  are not  common neighbours of $u$ and $v$.  So, as $|V(T)|=n-a+1$ and $U\cup V_\pi$ has been used, we are able to find
    a common neighbour of $u$ and $v$ in which to embed $x_j$.     This finishes the embedding of $T$.

It only remains to prove our claim that there is a permutation of $V(G)\setminus U$ such that~\eqref{A} and~\eqref{B} hold. We take
    a random permutation $\pi=(v_1,v_2,\ldots, v_{n'})$ of $V(G)\setminus U$, and show that with positive probability, it has both these properties. 
    We note that by~\eqref{degG'}, and since $V(G)\setminus U\subseteq V(G')$, we have
    $$\mathbb E\big[ |J_\pi|\big]\le \sum_{v\in V_\pi} \frac{|V(G)\setminus N(v)|}{|V_\pi|} \le |V_\pi|\cdot \frac {4\sqrt \delta k + a}{|V_\pi|}\le 5\sqrt\delta k.$$
    Hence by Markov's inequality  (see~\cite{mr02}), the probability that~\eqref{A} fails is at most $\frac 16$.

We will show that~\eqref{B} fails with probability less than $\frac 56$, which will finish the proof of our claim. 
    By definition of $H$, each vertex from  $H$ has  less than $a$ non-neighbours in $V(G')\setminus U$. Moreover, as $|U|\ge |S'|>a$, for any $v\in V(G')\setminus U$, the probability that $v\notin V_\pi$ is at most~$\frac{1}{50}$.
    So, setting $$s_\pi:=\sum_{u\in H\setminus (U\cup V_\pi)} |\{v \ |\ u \notin U,  uv\notin E(G')\}   \setminus V_\pi|,$$ we have that $$\mathbb E[s_\pi]\le \frac{|H\setminus U|a}{2500} \le \frac{ak}{2400},$$ 
    and by Markov's inequality  (see~\cite{mr02}), the probability that $s_\pi\ge \frac{ak}{1600}$ is at most~$\frac 23$. In particular, the probability that $H\setminus (U\cup V_\pi)$ has  more than $\frac k{500}$ vertices which each have at least $\frac{a}{3}$ non-neighbours in 
    $V(G')\setminus (U\cup V_\pi)$ is at most $\frac 23$. So, if we can show that 
    the probability that $|H\setminus (U\cup V_\pi)|<\frac k{500}+16\sqrt\delta k$ is at most $\frac 16$, we are done.
    
    For this, note that $|H\setminus U| 
     >\frac 4{25}k$ by~\eqref{H} and~\eqref{U}.
Also by~\eqref{U}, for each $v\in V(G')\setminus U$, the probability that $v\notin V_\pi$ is at least~$\frac{1}{60}$. It follows that the expectation of $|H\setminus (U\cup V_\pi)|$ is  at 
    least $\frac{k}{375}$. Applying the Chernoff bound (see~\cite{mr02}), 
    we deduce that the probability that  $|H\setminus (U\cup V_\pi)|<\frac{k}{400}$ is (much) less than $\frac 16$. As $\frac k{400}>\frac k{500}+16\sqrt\delta k$, we are done.

\newcommand{\etalchar}[1]{$^{#1}$}

\end{document}